\documentstyle[12pt]{article}
\newcommand{\const}{\mathop{\rm const}\limits}

\newcommand{\msupp}{\mathop{\rm msupp}\limits}

\newcommand{\vraisup}{\mathop{\rm vraisup}\limits}

\begin{document}
\begin{center}

{\bf Strichartz \ - \ type Inequalities for Parabolic and }\\

\vspace{3mm}

{\bf Schr\"odinger Equations in rearrangement invariant Spaces }\\

\vspace{3mm}

{\sc Ostrovsky E., Rogover E.}\\

\normalsize

\vspace{3mm}
{\it Department of Mathematics and Statistics, Bar \ - \ Ilan University,
59200, Ramat Gan, Israel.}\\
e-mail: \ galo@list.ru \\

\vspace{3mm}

{\it Department of Mathematics and Statistics, Bar \ - \ Ilan University,
59200, Ramat Gan, Israel.}\\
e - mail: \ rogovee@gmail.com \\

\vspace{3mm}

 {\bf Abstract.} \\

\end{center}

 In this paper we generalize the classical Strichartz
estimation for solutions of initial problem for linear parabolic and
Schr\"odinger PDE on  many popular classes {\it pairs} of rearrangement invariant(r.i.) spaces  and construct some examples in order to show the exactness of our estimations.\\

\vspace{3mm}


\vspace{3mm}


\vspace{3mm}

{\bf 1. Introduction. Notations. Statement of problem.}\\

\vspace{3mm}

{\bf Problem (P). } Let us consider the initial (Cauchy) problem for
the non \ - \ degenerate
linear parabolic  equation in the whole d \ - \ dimensional space $ R^d:$

$$
\frac{\partial u}{ \partial t} = \sum_{k =1}^d \sum_{j=1}^d a_{k,j}(t,x)
\frac{\partial^2 u}{\partial x_k \ \partial x_j }, \eqno(1.1)
$$
where $ u = u(t,x), \  t \in [0, \infty), \  x = \vec{x} = \{x_1, x_2,
\ldots, x_d \} $ be a $ d \ - \ $ dimensional vector:  $ x \in R^d $ with initial condition

$$
\lim_{t \to 0+} u(t,x) = f(x), \eqno(1.2)
$$
where the limit is understood in the $ L_p $ sense for some $ p \in [1,\infty].$
\par
 It is presumed that for the problem (P) (1.1) with  (1.2) are satisfied the classical conditions for existence and uniqueness, for instance:

$$
\lambda \sum_{k=1}^d \xi^2_k \le
\sum_{k=1}^d \sum_{j=1}^d a_{k,j}(t,x) \ \xi_k \xi_j \le \Lambda \sum_{k=1}^d \xi^2_k  \eqno(1.3)
$$
for some constants $ \lambda, \ \Lambda: \ 0 < \lambda \le \Lambda < \infty $
(the uniform ellipticity and boundeness condition);

$$
\max_k \max_j |a_{k,j}(t,x) \ - \  a_{k,j}(s,y)| \le M [|t \ - \ s| +
\sum_{l=1}^d |x_l \ - \ y_l|]^{\beta} \eqno(1.4)
$$
for some finite positive constants $ M, \beta; \ \beta \in (0,1] $ (the uniform H\"older condition). \par

 We denote as usually for arbitrary measurable (complex, in general case)
function $ g: R^d \to R $ and  for $  p = \const \in [1, \infty) \
p^/ \stackrel{def}{=} p/(p \ - \ 1), \ p = +\infty \ \Rightarrow p^/ = 1; $

$$
|g|_p = \left[\int_{R^d} |g(x)|^p \ dx \right]^{1/p}; \ L_p =
\{ g: \ |g|_p < \infty \},
$$

$$
L_{\infty} = \{g, \ \vraisup_{x \in R^d} |g(x)| \stackrel{def}{=}
|g|_{\infty} \} < \infty.
$$

 We will denote hereafter as $ C, \ C_m, \ C_j(d), \ C_l(d,a) $ some finite positive {\it non \ - \ essential} constants. \par

  It is well \ - \ known (\cite{Kapitanskii1}, \cite{Ginibre1} etc.)
 for the solution of problem (1.1) \ - \ (1.2) under conditions (1.3) \ - \ (1.4) the Strichartz  \ - \ Krylov estimations, which we want to reformulate in the convenient for us form. \par

 We will denote for the solution of (1.1) \ - \ (1.2)

$$
u(t,\cdot)  = T_t \ f(\cdot),
$$
where  $ \{ T_t \} $ is a semi \ - \ group of linear operators. \par

{\bf Lemma 1.1.} For all the values $ r > p, \ p \ge 1 $ and $ t > 0 $

$$
|u|_r \le C(a,d) \ |f|_p \ t^{\frac{d}{2} \left( \frac{1}{r} \ - \
\frac{1}{p} \right) }. \eqno(1.5)
$$

 The estimation (1.5) may be obtained as follows. Without loss of
generality we can and will assume the function $ x \to f(x) $ to be
non \ - \ negative  and non \ - \ trivial: $ |f|_p \in (0, \infty). $ \par

 The solution $ u = u(t,x) $ has a view:

$$
u(t,x) = C_1 \ t^{-d/2} \ \int_{R^d} f(y) \ G(t,x, y) \ dy, \eqno(1.6)
$$
where the positive
function $ G = G(t,x,y)$ is called Heat Potential (HP) and allows
the estimation: $ G(t,x,y) \le G_0(t,x \ - \ y); $

$$
G_0(t,z) = C_2 \ \exp \left( - C_3 \ t^{-1} \ ||z||^2 \right);  \
||z||^2 \stackrel{def}{=} \sum_{k=1}^d |z_k|^2.
$$

 So, we have:

$$
|u(t,x)| \le C_3 \ t^{-d/2} \ f*G_0, \eqno(1.7)
$$
where the convolution $ \ * \ $ is understood  over the variable $ x; \
t = \const > 0. $\par

 We obtain using the well \ - \ known Hardy \ - \ Littlewood \ - \ Young
inequality:

$$
|u|_r \le C_4 \ t^{-d/2} \ |f|_p  \ |G_0|_q, \  1 + 1/r = 1/p + 1/q,
$$

$$
 p,q,r \ge 1, \ r \ge p \ge 1.  \eqno(1.8)
$$

 It is easy to verify by the direct calculation that

$$
|G_o|_q \le C_5(d) \ t^{d/(2q)},
$$
therefore

$$
C_6 \ |u|_r \le t^{  d( \ - \ 0.5 + 0.5/q) } \ |f|_p \ =
t^{ 0.5 \cdot d( 1/r \ - \ 1/p) } \ |f|_p.
$$

{\bf Problem (S). } Let us consider also the initial (Cauchy) problem for
the (linear) Schr\"odinger's equation  without potential
("free particle") in the whole space $ R^d: $

$$
- i \frac{\partial v }{\partial t } =  0.5 \
\sum_{k=1}^d \frac{\partial^2 v}{\partial x_k^2 }
\stackrel{def}{=} 0.5 \ \Delta \ v, \  (i^2 = -1),
\eqno(1.9)
$$

$$
\lim_{t \to 0} u(t,x) = f(x) \eqno(1.10)
$$
again in the $ L_p $ sense for some $ p \in [1,\infty]. $ \par

 It is well \ - \ known (\cite{Strichartz1}, \cite{Strichartz2}, \cite{Bouclet1},
\cite{Kapitanskii1}, \cite{Bourgain1}, \cite{Chen1}, \cite{Killip1},
\cite{Rogers1}, \cite{Rogers2}, \cite{Rogers3} etc.)
for the solution $ v = v(t,x) $ of problem (S) (1.9) \ - \ (1.10) under
condition $ f(\cdot) \in L_q $ for some $ q \in [1,2] $ there exists, is
unique and satisfies the following assertion. \par

{\bf Lemma 1.2. } In the case $ p \ge 2 $ the following inequality is true:

$$
|v|_p \le C_7(d) \ |t|^{d(0.5 \ - \ 1/p^/) } \ |f|_{p^/}. \eqno(1.11)
$$

 The assertion (1.11) follows from the conservation law:
$$
\forall t > 0 \ \Rightarrow  \ |v|_2 \le C_8(d) \ |f|_2,
$$
explicit formula for $ v(\cdot,\cdot): t > 0 \ \Rightarrow $

$$
v(t,x) = C_8(d) \ t^{-d/2} \
\int_{R^d} \exp \left(0.5 \ i \ ||x \ - \ y||^2 /t \right) \ f(y) \ dy,
$$
from which it follows the inequality

$$
|v|_{\infty} \le C_8(d) \ t^{-d/2} \ |f|_1, \ t > 0,
$$
and from the famous interpolation theorem belonging to Riesz \ - \ Thorin.\par

 We will denote also for the solution $ v = v(t,x) $ of the problem (1.9) \ - \ (1.10)

$$
v(t,\cdot)  = U_t \ f(\cdot),
$$
where  $ \{ U_t \} $ is a group of linear operators. \par

 In the physical literature the operators  $ \{ T_t \} $ and $ \{ U_t \} $
are called often the {\it Propagation operators}.\par

 Note that the inequalities (1.5) and (1.11)  are non \ - \ trivial only in
the case  of sufficiently great values $ t: \ t >>1, $ or more exactly

$$
t \to \infty; \ t > 2.  \eqno(1.12)
$$

 Further we will assume the condition (1.12) to be satisfied for both considered problems.\par

\vspace{2mm}

{\bf Our goal is generalization of  the estimations (1.5) and (1.11)
on some popular classes of rearrangement invariant (r.i.) spaces, more exactly,
so \ - \ called moment rearrangement invariant (m.r.i.) spaces.}\par

\vspace{2mm}

{\it In detail. Parabolic case.}  The inequality (1.5)
 may be rewritten as follows. Let $ (X, ||\cdot||X) $ be any
  rearrangement invariant (r.i.) space on the set $ R^d; $ denote by
$ \phi(X, \delta) $  its fundamental function

  $$
  \phi(X,\delta) = \sup_{A, \mu(A) \le \delta} ||I(A)||X, \ I(A)= I(A,x) =
   1, x \in A,
  $$
  $ I(A) = I(A,x) = 0, \ x \notin A.$ \par

 We introduce also for two function r.i. spaces $ (X, \ ||\cdot||X ) $ and
$ (Y, \ ||\cdot||) $ defined over our set $ R^d \ $ and for arbitrary finite positive constants $ K_1,\ K_2 $ and the values $ t > 2 $
the so \ - \ called {\sc Strichartz Parabolic two \ - \ space functional, briefly: SP functional} between the spaces $ X $ and $ Y $ as

  $$
  W_{SP} (X,Y, K_1, K_2;t) \stackrel{def}{=} \sup_{f \in X, \ f \ne 0}
\left[ \frac{||T_t \ f||Y}{\phi(Y,K_1 t^{d/2}) } :
\frac{||f||X}{\phi(X,K_2 t^{d/2})} \right],
  $$
  $ W_{SP}(X,Y;t) \stackrel{def}{=} W_{SP}(X,Y,1,1; t).$ Then (1.5) is equivalent to the following inequality:

  $$
  r > p \ge 1 \Rightarrow  \sup_{t > 2} W_{SP}(L_p, L_r; t) < \infty. \eqno(1.13)
  $$

 \vspace{3mm}

{\bf Definition 1.} \par

\vspace{3mm}

{\sc By definition, the {\it pair} of r.i. spaces $ (X, \ ||\cdot||X ) $ and
 $ (Y, \ ||\cdot||) $ over $ R^d $ is said to be a {\it (strong) Strichartz
Parabolic pair}, write: $ (X,Y) \in SP, $ if the $ SP $ functional
$ W_{SP}(X,Y; t) $ between $ X $ and  $ Y $  is uniform on the variable
$ t, \ t > 2  $ finite:

  $$
  \sup_{t > 2} W_{SP}(X,Y;t) < \infty \eqno(1.14)
  $$
 and is called a {\it weak Strichartz Parabolic pair}, write $ (X,Y) \in wSP, $
 if for some non \ - \ trivial constants } $ K_1, K_2 $

  $$
 \sup_{t > 2} W_{SP}(X,Y, K_1, K_2;t) < \infty. \eqno(1.15)
  $$

\vspace{3mm}

{\it Schr\"odinger case.} \par

 We introduce also for two function r.i. spaces $ (X, \ ||\cdot||X ) $ and
$ (Y, \ ||\cdot||) $ defined over our set $ R^d \ $ and for arbitrary finite positive constant $ K $ and the values $ t > 2 $
the so-called {\sc Strichartz Schr\"odinger  two-space functional, briefly: SR
functional} between the spaces $ X $ and $ Y $ as

  $$
  V_{SR}(X,Y, K; t) \stackrel{def}{=} \sup_{f \in X, \ f \ne 0}
 \frac{ t^{-d/2} \ \cdot ||U_t \ f ||Y}{||f||X \ \cdot \phi(X, K t^{-d}) },
  $$
and define $ V_{SR}(X,Y; t) \stackrel{def}{=} V_{SR} (X,Y, 1; t). $ \par

\vspace{3mm}

{\bf Definition 2.} \par

\vspace{3mm}

{\sc By definition, the {\it pair} of r.i. spaces $ (X, \ ||\cdot||X ) $ and
 $ (Y, \ ||\cdot||) $ over $ R^d $ is said to be a {\it (strong) Strichartz
Schr\"odinger pair}, write: $ (X,Y) \in SR, $ if the $ SR $ functional
$ V_{SR}(X,Y; t) $ between $ X $ and  $ Y $  is uniform on the variable
$ t, \ t > 2 $ finite:

  $$
 \sup_{t > 2} V_{SR}(X,Y; t) < \infty \eqno(1.16)
  $$
 and is called a {\it weak Strichartz Schr\"odinger pair}, write $ (X,Y) \in wSR, $ if for some positive non \ - \ trivial constants } $ K $

  $$
   \sup_{t > 2} V_{SR}(X,Y, K; t) < \infty. \eqno(1.17)
  $$

  Roughly speaking, we will prove that the most of popular {\it pairs} of r.i. spaces are strong, or at last weak Strichartz pairs, Parabolic or
Schr\"odinger. \par

 The paper is organized as follows. In the next section we recall and describe
a {\it new class } of r.i. spaces, namely, so \ - \ called moment rearrangement
invariant spaces, briefly, m.r.i. spaces. In the section 3 we formulate and prove the main result of paper for m.r.i. spaces. \par
 In the section 4 we offer some examples of our results. In the section 5
we consider some low bounds for introduced functionals
in order to show the precision  of obtained estimations. \par
 In the last section 6 we describe some generalizations of results of the section 3. \par

\vspace{4mm}

{\bf 2. Auxiliary facts. Moment rearrangement invariant spaces.}\par

\vspace{3mm}

 The complete investigation of the theory of r.i. spaces see, e.g., in
\cite{Bennet1}, chapters 1,2; \cite{Krein1}, chapter 1. \par
 We recall here only that the Banach function space $ X $ equipped with the
norm $ || \cdot ||X $ over the set, e.g., $ R^d $ is called rearrangement
invariant (r.i.) space, if the norm in this space dependent only on the
distribution function of $ f: $

$$
||f||X = R(Q_f(\cdot)),
$$
where $ Q_f(\cdot) $ is  the distribution function for the (measurable)
function $ f: $

$$
Q_f(s) = m \{x, x \in R^d, \ |f(x)| > s \}; \ s \in (0, \infty)
$$
and $ R(\cdot) $ is some functional.\par

 For instance, many popular functional spaces: $ L_p $ spaces,
Orlicz, Lorentz, Marzinkiewitz spaces  are r.i. spaces.\par
  Let $ (X, \ ||\cdot||X) $ be the r.i. space, where $ X $ is linear subset on the space of all measurable function $ R^d \to R $  with norm $ ||\cdot||X. $
\par

\vspace{3mm}

 {\bf Definition 3.} \par

\vspace{3mm}

 {\sc We will say that the space $ X $ with the norm $ ||\cdot||X $ is {\it moment rearrangement invariant space,} briefly: m.r.i. space, or
$ X =(X, \ ||\cdot||X) \in m.r.i., $
 if there exist a real constants $ a, b; 1 \le a < b \le \infty, $ and some {\it rearrangement invariant norm } $ < \cdot > $ defined on the space of a real functions defined on the interval $ (a,b), $ non necessary to be finite on all the functions, such that

  $$
  \forall f \in X \ \Rightarrow || f ||X = < \ h(\cdot) \ >, \ h(p) = |f|_p. \eqno(2.1)
  $$

   We will say that the space $ X $ with the norm $ ||\cdot||X $ is {\it weak moment rearrangement space,} briefly, w.m.r.i. space, or $ X =(X, \ ||\cdot||X) \in w.m.r.i.,$ if there exist a constants $ a, b; 1 \le a < b \le \infty, $ and some {\it functional } $ F, $ defined on the space of a real functions defined on the interval $ (a,b), $ non necessary to be finite on all the functions, such that}
  $$
  \forall f \in X \ \Rightarrow || f ||X = F( \ h(\cdot) \ ), \ h(p) = |f|_p.  \eqno(2.2)
  $$

   Roughly speaking, the functional space $ X $ is called  m.r.i. space or
w.m.r.i. space, if the norm in this space dependent only on the some family
of $ L_p $ norms of considering function.\par

    We will write for considered w.m.r.i. and m.r.i. spaces $ (X, \
  ||\cdot||X) $

   $$
      (a, b) \stackrel{def}{=} \msupp(X), \eqno(2.3)
   $$
   ("moment support"; not necessary to be uniquely defined)
   and define for other such a space $ Y = (Y, \ ||\cdot||Y ) $ with
   $ (c,d) = \msupp(Y) $

   $$
   \msupp(Y) >> \msupp(X),
   $$
   or equally, $ \msupp(X) << \msupp(Y),$
    iff $ \max(a,b) \le \min(c,d). $ \par
     It is obvious that arbitrary m.r.i. space is r.i. space.\par

   There are many r.i. spaces satisfied the condition (2.2): exponential Orlicz's spaces, some Martzinkiewitz spaces, interpolation spaces (see
\cite{Astashkin1}, \cite{Jawerth1}, \cite{Davis1}, \cite{Steigenwalt1},
\cite{Ostrovsky5} etc. )\par
  In the article \cite{Lukomsky1}
 are introduced the so \ - \ called $ G(p,\alpha) $ spaces
  consisted on all the measurable function $ f: T \to R $  with finite norm
   $$
   ||f||_{p,\alpha} = \left[ \int_1^{\infty} \left(\frac{|f|_x}{x^{\alpha}}
    \right)^p \ \nu(dx) \right]^{1/p},
   $$
   where $ \nu $ is some Borelian measure.\par
    Astashkin in \cite{Astashkin2} proved that the space $ G(p,\alpha) $ in the    case  $ T = [0,1] $ and $ \nu = m, \ m $ is usually Lebesque measure,
    coincides with the
   Lorentz $ \Lambda_p( \log^{1-p \alpha}(2/s) ) $ space.  Therefore,
     both this spaces are m.r.i. spaces.\par

     Another examples. Recently (see \cite{Davis1}, \cite{Fiorenza1},
      \cite{Fiorenza2}, \cite{Ivaniec1}, \cite{Ivaniec2},
     \cite{Kozachenko1}, \cite{Ostrovsky1}, \cite{Ostrovsky2},   \cite{Ostrovsky3}, \cite{Ostrovsky4}, \cite{Ostrovsky5},\cite{Ostrovsky6} etc.)
     appear the so-called Grand Lebesque Spaces $ GLS = G(\psi) =
    G(\psi; a,b) $ spaces consisting on all the measurable functions
     $ f: R^d \to R $ with finite norms

     $$
     ||f||G(\psi) \stackrel{def}{=} \sup_{p \in (a,b)}
     \left[ |f|_p /\psi(p) \right]. \eqno(2.4)
     $$

      Here $ \psi(\cdot) $ is some positive continuous on the {\it open}      interval  $ (a,b) $ function such that

     $$
     \inf_{p \in (a,b)} \psi(p) > 0. \eqno(2.5)
     $$
      It is evident that $ G(\psi; a,b) $ is m.r.i. space and
   $ \msupp(G(\psi(a,b)) = (a,b). $\par

  We will write in this case $ \psi \in \Psi(a,b). $ \par

  This spaces are used, for example, in the theory of probability
   ( \cite{Kozachenko1}, \cite{Ledoux1},
     \cite{Ostrovsky1}, \cite{Ostrovsky2}, \cite{Ostrovsky3},
      \cite{Ostrovsky4}, \cite{Ostrovsky5},\cite{Ostrovsky6} etc.),
  theory of PDE (\cite{Fiorenza2},  \cite{Ivaniec2}),
 functional analysis (\cite{Astashkin1}, \cite{Astashkin2},
\cite{Davis1}, \cite{Fiorenza1}),  \cite{Jawerth1},
theory of Fourier series (\cite{Ostrovsky5}), theory of martingales
(\cite{Ostrovsky1}), \cite{Ostrovsky6}) etc.\par

 We can consider the classical Lebesgue spaces $ L_s, \ s \ge 1 $ as an "extremal" case of $ G(\psi) $ spaces, namely, define a function

$$
\psi_s(p) = 1, \ p = s, \ \psi_s(p) = + \infty, \ p \ne s.
$$

 If we define formally $ \infty/\infty = \infty, $ then
$$
|f|_s = ||f||G(\psi_s).
$$
 See in detail \cite{Ostrovsky2}, chapters 1,2. \par

  Let us consider as an example
now the (generalized) Zygmund's spaces $ L_p \ Log^r L, $ which may be defined as an Orlicz's spaces over some subset of the space
$ R^d $  with non-empty interior  and with $ N - $ Orlicz's function of a
 view

    $$
    \Phi(u) = |u|^p \ \log^r( C + |u|), \ p \ge 1, \ r \ne 0.
    $$

     \vspace{3mm}

     {\bf Lemma 2.1.} \\

    \vspace{2mm}

     {\bf 1.} All the spaces $ L_p \ Log^r L $ over real line with measure $ m $
     with condition $ \ r \ne 0 $ are not m.r.i. spaces.\\
     {\bf 2.} If $ r $ is positive and integer, then the spaces $ L_p Log^r L $    are w.m.r.i. space.\\

\vspace{2mm}

  {\bf Proof.} {\bf 1.} It is sufficient to consider the case $ d = 1 $   with the classical Lebesgue 
measure $ m $ and the case $ p > 1. $ \par
      There exists a function $ f_0 = f_0(x) $ belonging to the space
    $ L_p \ Log^r L,  $ for example, for which 
  $$
  \int_T |f_0|^p \ \log^r(C + |f_0|) \ dx < \infty,
  $$
  but such that for all sufficiently small values $ \epsilon > 0 $

  $$
  \int_T |f_0|^{p \pm \epsilon} \ dx = \infty
  $$
 in the case $ p > 1 $ and

  $$
  \int_T |f_0|^{p + \epsilon} \ dx = \infty
  $$
   in the case $ p = 1. $ \par

 Therefore, the interval $ (a,b) $ in the definition of m.r.i. spaces does not exists.\par
  The affirmation {\bf 2} it follows from the formula 

  $$
  |f|^p \ [\log|f|]^k = d^k |f|^p /dp^k, \ k = 1,2,\ldots.
  $$

\vspace{3mm}

  {\bf Lemma 2.2} There exists an r.i. space without the w.m.r.i. property.\par

\vspace{3mm}

  {\bf Proof.} On the interval $ T = [0,1] $ with usual Lebesque measure $ m $ there exists a function $ f $ with standard normal (Gaussian) distribution.
This implies, for example, that

$$
\int_T \exp(p f(x)) \ dx = \exp \left(0.5 \ p^2 \right), \ p \in R.
$$

There exist a functions $ g: R \to R $ such that the function $ h(x) = g(f(x)) $ which {\it distribution} can not be uniquely
defined by means of all positive moments, for instance, $ h(x) = g(f(x)) = [f(x)]^3 $ or $ g(x) = \exp(f(x)). $ \par
   Let us consider a two such a functions $ f_1 $ and $ f_2 $ with {\it different} distributions, but with at the same moments, for example:

   $$
   \int_T |f_1|^p \ dx = \int_T |f_2|^p \ dx = \int_T [\exp(f)]^p \ dx = \exp(p^2/2), \ p \in R.
   $$

 We choose the (quasi) \ - \ concave positive strictly increasing continuous function $ \theta(\cdot), \ \theta(0+) = 0, $ for which

 $$
 \int_0^{\infty} \theta( m \{x: \ |f_1(x)| > \lambda \} ) \ d\lambda = \infty,
 $$
 but
 $$
 \int_0^{\infty} \theta( m \{ x: \ |f_2(x)| > \lambda \} ) \ d\lambda < \infty.
 $$
 The Lorentz r.i. space $ \Lambda(T, \theta) $ over $ T = [0,1] $ with the function
 $ \theta(\cdot) $ and the classical norm (see \cite{Bennet1}, chapter 2, section 2)

  $$
  ||f||L(T,\theta) = \int_0^{\infty} \theta( m \{x: \ |f(x)| > \lambda \} ) \ d \lambda
  $$
 is not w.m.r.i. space.\par

\vspace{5mm}

  {\bf 3. Main result. Strichartz inequalities for the pairs of m.r.i. spaces.}\\

\vspace{2mm}

{\bf Theorem 3.1.} {\it Let $ (X, ||\cdot||X) $ be any m.r.i. space over the space $ R^d $
with moment support $ \msupp(X) =(a_1,b_1) $ relatively the auxiliary norm $ < \cdot >, $ and let $ (Y,||\cdot||Y ) $ be
another m.r.i. space over at the same set $ R^d $ relatively the
second auxiliary norm $ << \cdot >> $ and with $ \msupp(Y) = (a_2,b_2), $ where
$ (a_1,b_1) << (a_2,b_2). $ \par
  Then the pair of m.r.i. spaces $ (X, ||\cdot||X) $ and $ (Y, ||\cdot||Y) $
is the (strong) Parabolic Strichartz pair: }

$$
\sup_{t > 2} \ W_{SP}(X,Y;t) = C_{SP}(X,Y) < \infty. \eqno(3.1)
$$

\vspace{2mm}

{\bf Note } that the restriction $ (a_1,b_1) << (a_2,b_2) $
is not loss of generality.\par

\vspace{3mm}

{\bf Theorem 3.2.}
{\it Let $ (X, ||\cdot||X) $ be any m.r.i. space over the space $ R^d $
with moment support $ \msupp(X) =(a_1,b_1) $ relatively the auxiliary norm $ < \cdot >, $ and let $ (Y,||\cdot||Y ) $ be
another m.r.i. space over at the same set $ R^d $ relatively the
second auxiliary norm $ << \cdot >> $ and with $ \msupp(Y) = (a_2,b_2), $ where
$ (a_1,b_1) << (a_2,b_2), \ a_2 > 2. $ \par
  Then the pair of m.r.i. spaces $ (X, ||\cdot||X) $ and $ (Y, ||\cdot||Y) $
is the (strong) Schr\"odinger Strichartz pair: }

$$
\sup_{t > 2} \ V_{SR}(X,Y;t) = C_{SR}(X,Y) < \infty. \eqno(3.2)
$$

{\bf Proofs. Theorem 3.1.} It follows from the inequalities (1.5) for the
values  $  p \in (a_1, b_1) $ and $ r \in (a_2, b_2) $ correspondingly:

$$
|T_t \ f|_r \ t^{d/(2p)} \le C(a,d) \  |f|_p \ t^{d/(2r)}. \eqno(3.3)
$$

Tacking into account the monotonicity  of the norm $ < \cdot > $ and equality

$$
      \phi(X,t^{d/2}) = < z(\cdot) >,
$$
where  for all admissible values $ r $

$$
z( r ) = t^{d/(2r)},
$$
we get from (3.3) tacking the norm $ < \cdot >: $

$$
|T_t \ f|_r \ \cdot \phi(X,t^{d/2}) \le C(a,d) \ ||f||X \ \cdot t^{d/(2r)}. \eqno(3.4)
$$

 Tacking analogously
from the bide \ - \ side  of inequality (3.4) the norm $ << \cdot >>,$
we conclude

$$
||T_t \ f||Y \cdot  \phi(X,t^{d/2}) \le C(a,d) \ ||f||X \cdot \phi(Y,t^{d/2}),
$$
which completes the proof of Theorem 3.1. \par

\vspace{3mm}

{\bf Proof of Theorem 3.2} is analogously. We use the assertion (1.11)
of the Lemma 1.2:

$$
|U_t \ f|_p \ t^{d/(2p^/)} \le C(d) \ t^{d/2} \ |f|_{p^/}.
$$

Tacking the norm $ < \cdot >, $ we obtain:

$$
|U_t \ f|_p \ \cdot \phi(X, t^{d/2}) \le C(d) \ t^{d/2} \ \cdot || f ||X.
$$

Tacking the norm $ << \cdot >>, $ we obtain:

$$
||U_t \ f||Y \  \cdot \phi(X, t^{d/2}) \le C(d) \ t^{d/2} \ \cdot || f ||X.
$$

This completes the proof of Theorem 3.2. \par

Note now as a particular case the case when $ X = G(\psi), \ \psi \in
\Psi(a_1, b_1); \ Y = g(\nu), \ \nu \in \Psi(a_2,b_2), \ b_1 < a_2: $

$$
\sup_{f \in G(\psi), \ f \ne 0} \ \sup_{t > 2}
\left[ \frac{||T_t \ f||G(\nu)}{\phi(G(\nu), t^{d/2})} :
\frac{||f||G(\psi)}{\phi(G(\psi), t^{d/2})} \right] = C_1(\psi, \nu) < \infty \eqno(3.5)
$$
and

$$
\sup_{f \in G(\psi), \ f \ne 0} \ \sup_{t > 2}
\frac{ t^{-d/2} \ ||U_t \ f ||G(\nu)}{||f||G(\psi) \cdot \phi(G(\psi), t^{-d}) } = C_2(\psi, \nu) < \infty.  \eqno(3.6)
$$

\vspace{3mm}

{\bf 4. Examples. } \par

\vspace{4mm}

 We consider now a very important for applications examples of $ G(\psi) $
spaces. Let $ a = const \ge 1,
b = const \in (a, \infty]; \alpha, \beta = const. $ Assume also that at
$ b < \infty \ \min(\alpha,\beta) \ge 0 $ and denote by $ h $ the (unique)
root of equation
$$
(h-a)^{\alpha} = (b-h)^{\beta}, \ a < h < b;
  \ \zeta(p) = \zeta(a,b; \alpha,\beta; p) =
$$
$$
(p-a)^{\alpha}, \ p \in (a,h);
\ \zeta(a,b; \alpha,\beta;p) = (b-p)^{\beta}, \ p \in [h,b);
$$
and in the case $ b = \infty $ assume that
 $ \alpha \ge 0, \beta < 0; $ denote
by $ h $ the (unique) root of equation
 $ (h-a)^{ \alpha} = h^{ \beta}, h > a; $ define in this case

$$
\zeta(p) = \zeta(a,b;\alpha,\beta;p) = (p-a)^{\alpha}, \ p
\in (a,h); \ p \ge h \ \Rightarrow \zeta(p) = p^{\beta}.
$$

   Here and further
 $ p \in (a,b) \ \Rightarrow \psi(p) \asymp \nu(p) $ denotes that

$$
0 < \inf_{p \in (a,b)} \psi(p)/\nu(p) \le \sup_{p \in (a,b)} \psi(p)/\nu(p)
< \infty.
$$

 The space $ G = G_{R^d}(a,b;\alpha,\beta)=
G(a,b; \alpha,\beta) $ consists by definition on all the measurable functions
$ f: T \to R $ with finite norm:
$$
||f|| G(a,b; \alpha,\beta)= \sup_{p \in (a,b)}
\left[ |f|_p \cdot \zeta(p) \right]. \eqno(4.1)
$$

 On the other words, $ G(a,b; \alpha,\beta)$ is the $ G(\psi; a,b) $ space with
$ \psi(p) = 1/\zeta(p). $ \par
 These spaces was introduced in \cite{Kozachenko1}, \cite{Ostrovsky2},
\cite{Ostrovsky5}; and in the two last articles was also calculated its fundamental functions. \par
 We rewrite here only the asymptotical expression for
$ \phi(G(a,b; \alpha,\beta) \ \delta) $ for two cases: $ \delta \to 0+ $ and
$ \delta \to \infty. $ \par

{\bf 1}. As $ \delta \to 0+: $

$$
\phi(G(a,b; \alpha,\beta),\delta) \sim (\beta b^2/e)^{\beta} \cdot
\delta^{1/b} |\log \delta|^{-\beta}, \eqno(4.2)
$$

$ 1 \le a < b < \in \infty, \alpha, \beta \ge 0; $

$$
\phi(G(c,\infty; \alpha, \ - \ \beta), \delta) \sim (\beta)^{|\beta|}
|\log \delta|^{-|\beta|}, \eqno(4.3)
$$

$ 1 \le c, \ \alpha \ge 0, \beta > 0; $

\vspace{3mm}

{\bf 2}. As $ \delta \to \infty: $

$$
\phi(G(a,\infty; \alpha, \ - \ \beta),\delta) \sim ( a^2 \alpha/e)^{\alpha}
\ \delta^{1/a} (\log \delta)^{-a}, \eqno(4.4)
$$

$ 1 \le a < b \le  \infty. $ \par

\vspace{3mm}

 We choose in this pilcrow

$$
X = G(a_1, b_1; \alpha_1, \beta_1), \ Y = G(a_2, b_2; \alpha_2, \beta_2),
$$
where $ 1 \le a_1 < b_1 < a_2 < b_2 \le \infty. $ \par

\vspace{3mm}

{\bf Parabolic example.}\par

\vspace{2mm}

 We obtain using the theorems 3.1  for the values $ t > 2:  $

$$
 \sup_{f: f \ne 0, f \in X} \ \left[
\frac{ ||T_t \ f||G(Y) }{ ||f||G(X) } \right] \le
$$

$$
C_1(d; a_1,b_1,a_2,b_2; \alpha_1,\alpha_2,\beta_1,\beta_2) \
t^{ - \frac{d}{2 } \left( \frac{1}{a_1 } \ - \ \frac{1 }{a_2 } \right) } \cdot
(\log t)^{\alpha_2 \ - \ \alpha_1}. \eqno(4.5)
$$

\vspace{3mm}

{\bf Schr\"odinger example.}\par

\vspace{2mm}

We consider again  the case when

$$
X = G(a_1, b_1; \alpha_1, \beta_1), \ Y = G(a_2, b_2; \alpha_2, \beta_2),
$$
where $ 1 \le a_1 < b_1 < a_2 < b_2 \le \infty, $ but assume in addition
$ b_1 \le 2, a_2 \ge 2 $ (the cases $ b_1 > 2 $ or $ a_2 < 2 $ are trivial).
\par

 We obtain using the theorems 3.2  for at the same values $ t > 2:  $

$$
 \sup_{f: f \ne 0, f \in X} \ \left[
\frac{ ||U_t \ f||G(Y) }{ ||f||G(X) } \right] \le
$$

$$
C_2(d; a_1,b_1,a_2,b_2; \alpha_1,\alpha_2,\beta_1,\beta_2) \
t^{\frac{d}{2} \ - \ \frac{d}{b_1^/ }} \  (\log t)^{-\beta_1}. \eqno(4.6)
$$

\vspace{3mm}

{\bf 5. Low bounds. } \par

\vspace{4mm}

 In this subsections we will construct some examples in order to illustrate the exactness of result of section 3, for example, the exactness of inequalities
(4.5) and (4.6). \par

{\bf Theorem 5.1.} {\it Let } $ X = L_1(R^d) $ {\it and }
$ Y = G(\nu), \ \nu $ {\it is arbitrary function from the space}
$ \Psi(a,b): \nu(\cdot) \in \Psi(a,b), a > 1,  a < b \le \infty $
{\it be two examples of r.i.}  {\it spaces. We assert that}

\vspace{3mm}

$$
\underline{\lim}_{t \to \infty} \sup_{f: f \ne 0, f \in X} W_{SP}(X, G(\nu);t ) = C^{(1)}_{SP} > 0. \eqno(5.1)
$$

\vspace{2mm}

{\bf Theorem 5.2.} {\it Let} $  X = L_1(R^d) $ {\it and }
$ Y = L_{\infty}(R^d) $ {\it be two examples of } $ G(\Psi) $ {\it spaces. We assert that}

$$
\underline{\lim}_{t \to \infty} \sup_{f: f \ne 0, f \in G(\psi)}
V_{SR}(X, Y ;t) = C^{(2)}_{SR} > 0. \eqno(5.2)
$$

\vspace{2mm}

{\bf Proof of theorem 5.1.} \par

\vspace{2mm}

{\bf 1.} It is sufficient to consider here in the problem (P) only the
case if equation (1.9) has a view
$$
 \frac{\partial u }{\partial t } =  0.5 \
\sum_{k=1}^d \frac{\partial^2 u}{\partial x_k^2 }.
\eqno(5.1)
$$

{\bf 2.} Let us consider the following function (multidimensional normal density)  for the values $ x = \vec{x} \in R^d, \ \sigma = \const > 0: $

$$
g(x) = g_{\sigma}(x) = (2 \pi \ \sigma)^{-d/2} \
\exp \left(-||x||^2/(2 \sigma^2) \right). \eqno(5.2)
$$

 We get after direct calculation  for $ q \in [1,\infty]: $

$$
|g|_q \asymp C(d) \ | \sigma|^{-d \left( \frac{1}{2} \ - \ \frac{1}{q} \right)}.
\eqno(5.3)
$$

{\bf 3.} Let $ \nu \in G(\Psi), \ \Psi \in (a,b) $ and $ f(x) = g(x). $
We have:

$$
||f||X = |g|_1 = 1; \ \phi(X; \delta) = \delta, \delta > 0;
$$
therefore

$$
\phi \left(X; t^{d/2} \right) = t^{d/2}.
$$

{\bf 4.} As long as the solution of the equation (5.1) has a view:

$$
u = u(t,x) = g_{1 + t}(x),
$$

we have for sufficiently greatest values $ t, \ t > 2 $  and $ r \in (a,b): $

$$
|u|_r \asymp t^{\frac{d}{2}( \frac{1}{r} \ - \ 1) };
$$

$$
||u||G(\nu) \asymp \sup_{r \in (a,b)} \left[t^{\frac{d}{2}
( \frac{1}{r} \ - \ 1 ) } \right]/\nu(r) =
t^{-d/2} \ \phi \left(G(\nu); t^{-d/2} \right).
$$

{\bf 5.} Substituting into the expression for $ W_{SP}(X, G(\nu);t ), $ we
conclude that for $ t \ge 2 $ and $ f = g $

$$
W_{SP}(X, G(\nu);t ) \ge C \ \left[ \frac{||T_t \ g||Y}{\phi(Y, \ t^{d/2})} :
\frac{||g||X}{\phi(X, \ t^{d/2})} \right] \ge
$$

$$
C \ \left[ \frac{ t^{-d/2} \ \phi \left(G(\nu); t^{-d/2} \right)}
{\phi(Y, t^{d/2}) } : \frac{1}{t^{d/2}} \right] =
$$

$$
C \ \left[  t^{-d/2} : \frac{1}{t^{d/2}} \right] = C > 0.
$$

{\bf Proof of theorem 5.2 .} We choose as a function $ f(x) =
u(0,x) $ again the function $ f(x) = g(x) $ and obtain:

$$
u = u(t,x) = g_{1 + it}(x).
$$

Note that the formula (5.3) remains true for the complex values $ \sigma,
|\sigma| \ge 2. $ \par
 We have:

$$
|g|_1 = 1, \ |g|_{\infty} \asymp C,
$$

$$
|U_t \ g|_1 \asymp C, \ |U_t \ g|_{\infty} \asymp t^{-d/2},
$$
and at $ r \in (1,\infty) $

$$
| U_t \ g|_r \asymp
t^{-d \left( \frac{1}{2} \ - \ \frac{1}{r} \right)}. \eqno(5.4)
$$

Therefore  at $ t > 2 $

$$
V \left(L_1, L_{\infty}; t \right) \ge C
\frac{ t^{-d/2} \cdot t^{-d/2}}{ t^{-d}} = C > 0,
$$
QED.\\

\vspace{3mm}

{\bf 6. Concluding remarks.} \par

\vspace{3mm}

{\bf A. Mix estimations.}\\

\vspace{3mm}

Let for $ T > 0 \ S = S_T = (0, T) $ and $ \theta = \theta(\cdot) = \theta(t),
 \ t \in S_T \in \Psi(A,B), \ 1 \le A < B \le \infty. $ We denote the norm on the space consisting on all the measurable functions $ h: S_T \to R $
$ G(\theta) $ as $ ||| \ h \ |||G(\theta) $ and introduce the so \ - \ called {\it mix norm}  as

$$
||| \ || T_t \ f || \ |||G(Y) \times G(\theta_T) \stackrel{def}{=}
||| \ ||u(\cdot, \cdot)||G(Y) \ |||G(\theta_T).
$$

\vspace{2mm}

{\bf Theorem 6.1.}  {\it Let} $ X, Y $ {\it be m.r.i. spaces such that} $ \msupp(X) << \msupp(Y). $ {\it It follows from the theorem 3.1 that}

$$
||| \ || u(\cdot, \cdot)|| \ |||G(Y) \times G(\theta_T) \le C \ ||f||X \ \cdot
 ||| \ \frac{\phi(Y, t^{d/2}) }{\phi(X, t^{d/2}) } \ ||| G(\theta_T), \eqno(6.1)
$$
if  of course the last norm $ ||| \ \cdot \ ||| $ is finite.\\

\vspace{3mm}

{\bf B. Generalizations. }\\

\vspace{3mm}

 Let us consider some generalization of Schr\"odinger equation of a view:

 $$
 \partial_t u + (- \Delta)^{\alpha/2}u = 0, \ 0 < \alpha \le 2. \eqno(6.3)
 $$

$$
u(0,x) = f(x),
$$
the so \ - \ called {\it dispersive equation}, {\it non \ - \ local
diffusion equation} or {\it model of Keller \ - \ Segel},
see \cite{Biler1}, \cite{Fino1}. \par

 We denote the (unique) solution of (6.3) as $ u = S_{\alpha}(t) \ f. $ \par
In the article \cite{Biler1} is proved the estimation for $ u(t,x) $ of a view:

$$
|S_{\alpha}(t) \ f|_r \le C t^{\frac{d }{\alpha }( \frac{1}{r} \ - \ \frac{1}{p})} \ |f|_p,
\eqno(6.4)
$$
$ 1 \le p \le r \le \infty,$  and

$$
|\Delta S_{\alpha}(t) \ f|_r \le C t^{\frac{d }{\alpha }( \frac{1}{r} \ -
\ \frac{1}{p}) \ - \ \frac{1}{\alpha}} \ |f|_p. \eqno(6.5)
$$

We conclude repeating the proof of theorem 3.1 and using the inequality 6.4: \par

{\bf Theorem 6.2} {\it We have under the condition of theorem 3.1}

$$
\sup_{f \in X, \ f \ne 0} \ \sup_{t > 2}
\left[ \frac{||S_{\alpha}(t) \ f||Y}{\phi(Y, t^{d/\alpha}) } :
\frac{||f||X}{\phi(X, t^{d/\alpha})} \right] = C_1(\alpha, X,Y) < \infty.
\eqno(6.6)
$$

 The inverse assertion to the theorem 6.2 is also true in the following 
sense: \\
 
{\bf Theorem 6.3.} {\it Let } $ \alpha = 2, \ X = L_1(R^d) $ {\it and }
$ Y = G(\nu), \ \nu $ {\it is arbitrary function from the space}
$ \Psi(a,b): \nu(\cdot) \in \Psi(a,b), a > 1,  a < b \le \infty $
{\it be two examples of r.i. spaces. We assert that}

\vspace{3mm}

$$
\underline{\lim}_{t \to \infty} \sup_{f: f \ne 0, f \in X}
\left[ \frac{||S_{2}(t) \ f||Y}{\phi(Y, t^{d/2}) } :
\frac{||f||X}{\phi(X, t^{d/2})} \right] = C_{(2),SR}(X,Y) > 0. \eqno(6.7)
$$

 The {\bf proof} used the inequality (5.4) and is completely alike to the
proof f the theorem 5.2. For instance, we can choose instead the function
$ f $ the function $ f(x) = g(x) $ etc. \par

\vspace{3mm}

{\bf C. Derivatives.}\\

\vspace{3mm}

{\bf Theorem 6.4} We have under the condition of theorem 3.1 using the estimation (6.5)

$$
\sup_{f \in X, \ f \ne 0} \ \sup_{t > 2} \ t^{1/\alpha} \cdot
\left[ \frac{||\Delta S_{\alpha}(t) \ f||Y}{\phi(Y, t^{d/\alpha}) } :
\frac{||f||X}{\phi(X, t^{d/\alpha})} \right] = C_3(\alpha, X,Y) < \infty. \eqno(6.8)
$$

 This  estimations (6.8) is exact as in the theorem 6.3.\\

\vspace{4mm}





 \vspace{4mm}


\begin{center}

AUTHORS \\

{\bf Ostrovsky E.}\\
\vspace{3mm}

 Address: Ostrovsky E., ISRAEL, 84105, Ramat Gan, Bar \ - \ Ilan University.\\

\vspace{3mm}

e - mail: {\bf Galo@list.ru}\\

\vspace{4mm}

{\bf Rogover E.}\\

\vspace{3mm}

 Address: Rogover E., ISRAEL, 84105, Ramat Gan, Bar \ - \ Ilan University. \\

\vspace{3mm}
e - mail: \ rogovee@gmail.com \\

\end{center}

\end{document}